\documentclass[11pt]{amsart}
\usepackage{amsfonts,amssymb}
\input xy
\xyoption{all}
\begin{document}

\newcommand{\A}{{\mathcal A}}
\newcommand{\C}{{\mathbb C}}
\newcommand{\comb}{{\rm HL}}
\newcommand{\gl}{{\rm Gl}}
\newcommand{\hh}{{\mathbb H}}
\newcommand{\link}{{\rm link}}
\newcommand{\oh}{{\rm O}}
\newcommand{\La}{{\langle}}
\newcommand{\Ra}{{\rangle}}
\newcommand{\pl}{{\rm PL}}
\newcommand{\pgl}{{\rm PGl}}
\newcommand{\T}{{\mathbb T}}
\newcommand{\tp}{{\rm Top}}
\newcommand{\wTheta}{{\widetilde{\Theta}}}
\newcommand{\Z}{{\mathbb Z}}
\newcommand{\Ell}{{\mathbb {L}_\pm}}
\newcommand{\tI}{{\tilde{I}}}

\title{String orientations of simplicial homology manifolds}
\author{Jack Morava}
\address{Department of Mathematics, Johns Hopkins University,
Baltimore, Maryland 21218}
\email{jack@math.jhu.edu}
\thanks{The author was supported in part by AIM and the NSF}
\date {31 Oct 2009}
\begin{abstract}{Simplicial homology manifolds are proposed as 
an interesting class of geometric objects, more general than 
topological manifolds but still quite tractable, in which questions 
about the microstructure of space-time can be naturally formulated. 
Their string orientations are classified by $H^3$ with coefficients in 
an extension of the usual group of D-brane charges, by cobordism
classes of homology three-spheres with trivial Rokhlin invariant.}
\end{abstract}

\maketitle

\noindent
[This note grew out of the May 08 AIM workshop on algebraic topology and
physics. Thanks to Hisham Sati for organizing that meeting, and to the participants for enduring 
multiple fragmentary presentations of these ideas. I owe M. Ando, D. Christensen, M. Furuta,
N. Kitchloo, A. Ranicki, R. Stern, L. Taylor and B. Williams particular thanks for help
(without their bearing any responsibilities for remaining howlers).] 
\medskip
                                             
\section{Notation and background}

\noindent
{\bf 1.1} In homotopy-theoretic terms, a {\bf string structure} on a smooth manifold $M$ is a lift 
\[
\xymatrix{
{} & B\oh \La 8 \Ra \ar[d] \\
M \ar@{.>}[ur] \ar[r] & B\oh }
\]                                                
of the map classifying the stable tangent bundle of $M$; where the 7-connected cover $B\oh \La
8 \Ra$ of $B\oh$ is the fiber 
\[
B\oh \La 8 \Ra \to B\oh \to (B\oh)_{(7)} 
\]
of the map to its Postnikov approximation [27] having homotopy groups concentrated in degrees
seven and below. Alternately:
\[
(B \oh) \La 8 \Ra = B(\oh \La 7 \Ra)
\]
is the classifying space for the topological group 
\[
\oh \La 7 \Ra \to \oh \to \oh_{(6)}
\]
whose homotopy groups agree with those of O in degrees greater than or equal to seven. \medskip

\noindent
A string structure can thus be interpreted as a `reduction' of the structure group of the stable
tangent bundle of $M$, from O to $\oh \La 7 \Ra$, just as a spin structure is similarly a
`reduction' of that structure group to
\[
{\rm Spin} = \oh \La 2 \Ra \to \oh \to \Z_2 \times H(\Z_2,1) \;.
\]
An oriented manifold $M$ admits a spin structure iff the Stiefel-Whitney class $w_2 = 0$; similarly, 
the existence of a string structure requires the vanishing of a version $p_1/2$ of the Pontrjagin 
class defined  for spin manifolds. The set of spin structures on $M$ admits a transitive free action of
$H^1(M,\Z_2)$, and by essentially the same argument the set of string structures is an
$H^3(X,\Z)$-torsor. The associated twisted $K$-groups of $X$ are natural repositories [7] for Ramond-Ramond 
$D$-brane charges. \medskip

\noindent
The interest in string orientations comes from from quantum field theory, where they were
recognized as necessary to define an analog of the Dirac operator on the space $LM$ of free loops
on $M$, but mathematical interest [1] in their properties goes back to the early 70's [16]. There is
a parallel interest in the representation of three-dimensional cohomology classes in local
geometric terms, analogous [6] to the description of two-dimensional classes by complex line
bundles. \bigskip

\noindent
{\bf 1.2} If $G$ is a connected, simply-connected simple Lie group, then 
\[
\pi_3(G) = \Z
\]
by a classical theorem of Bott. A theorem of Kuiper says that the group $\gl(\hh)$ of invertible
bounded operators on an (infinite-dimensional) Hilbert space is contractible; the projective
general linear group
\[
\pgl(\hh) = \gl(\hh)/\C^\times \sim B\T 
\]
is therefore an Eilenberg-MacLane space of type $H(\Z,2)$, from which it follows that the
classifying space for $\pgl(\hh)$-bundles is an EM space of type $H(\Z,3)$. The associated
bundle
\[
1 \to \pgl(\hh) \to G \La 4 \Ra \to G \to 1
\]
is an extension of topological groups [28]. The following construction is due to Kitchloo 
[18, appendix]: \medskip 

\noindent
A level one projective representation of the loop group of $G$ on $\hh$ defines a 
homomorphism to $\pgl(\hh)$ which pulls $\gl(\hh)$ back to a version of the universal
central extension of $LG$. Let $\A$ denote the topological Tits building of $LG$, modeled by
the contractible space of connections on a trivial $G$-bundle over a circle; the pointed loops
act freely on it, and the holonomy map makes it a principal $\Omega G$ bundle.
The diagonal action of $LG$ on
\[ 
\A \times_{\Omega G} \pgl(\hh) := G\La 4 \Ra
\]
factors through an action of $G$ on $G \La 4 \Ra$ lifting the action of $G$ on itself by conjugation.
\medskip

\noindent
Note that because $\pgl(\hh)$ is the group of automorphisms of $\gl(\hh)$, the cohomology
group
\[
H^1(X,\pgl(\hh)) \cong H^1(X,H(\Z,2)) \cong H^3(X,\Z)
\]
can be interpreted as a Brauer group of bundles of $C^*$-algebras (up to Morita equivalence)
over $X$. A refinement of this argument [21] represents the Brauer group of bundles of {\bf
graded} $C^*$-algebras over $X$ by the generalized Eilenberg-MacLane space
\[
H(\Z_2,1) \times H(\Z,3) \;,
\]
which suggests interpreting $\oh \La 7 \Ra = \oh \La 4 \Ra$ in terms of a bundle of
such algebras over (S)O.

\section{Simplicial homology manifolds}

\noindent
{\bf 2.1} I also want to thank Kitchloo for observing that in low dimensions, the map 
\[
\pl/\oh \to B\oh \to B\pl
\]
is almost an equivalence: the homotopy groups of the fiber are the Kervaire-Milnor 
groups of differentiable structures on spheres, which (aside from the smooth 4D Poincar\'e conjecture 
\dots) are trivial below dimension seven. This implies that a string structure on a smooth 
manifold is the same as a smoothing of a topological manifold endowed with a lift
\[
\xymatrix{
{} & B\tp \La 8 \Ra \ar[d] \\
M \ar@{.>}[ur] \ar[r] & B\tp }
\]                                                               
of the map classifying its tangent topological block bundle. The theorem
\[
\tp/\pl \sim H(\Z_2,3) 
\]
of Kirby and Siebenman [24] seems also to point in this direction. \bigskip

\noindent
{\bf 2.2} The classification of string structures on geometric objects more general than
smooth manifolds is accessible nowadays, thanks to many researcher-years of deep work
related to the Hauptvermutung, suggesting that questions like `Who ordered the differentiable
structure' may not be out of reach. I will summarize some background from Ranicki's elegant 
account [22], but in some cases I'll use terminology from [15]: \medskip

\noindent
{\sc Definition:} A space $X$ is a $d$-dimensional homology manifold iff for any $x \in X$,  
\[
H_*(X,X - \{x\};\Z) \cong H_*(S^{d-1};\Z) \;;
\]
but a {\bf simplicial} homology manifold is a simplicial complex
$K$ such that for any $k$-simplex $\sigma \in K$,
\[
H_*(\link_K(\sigma);\Z) \cong H_*(S^{d-k-1};\Z) \;.
\]
The polyhedron $|K|$ of $K$ is a homology manifold iff $K$ is a simplicial homology manifold. 
A manifold homology resolution $f : M \to X$ of a space $X$ is a compact topological manifold 
$M$ together with a surjective map $f$ with acyclic point inverses. \medskip

\noindent
The element $\varkappa_k(K) \in H^k(|K|,\Theta_{k-1})$ Poincar\'e dual to the cycle
\[
\sum_{|\sigma| = d-k} [\link_K(\sigma)] \cdot \sigma \in H_{d-k}(|K|,\Theta_{k-1}) 
\]
[29 p. 63-65] with coefficients in the group of simplicial homology spheres (up to
cobordism through PL homology cylinders) is trivial unless $k=4$: for $\Theta := \Theta_3$
is the only nontrivial coefficient group. [It is {\bf not} finitely generated [12, 15, 25].]
There is a block bundle theory [14, 15] for simplicial homology manifolds, resulting in a fibration 
\[
B\pl \to B\comb \to H(\Theta,4)
\]
of classifying spaces. \bigskip

\noindent
{\bf Theorem} [cf. [22 \S 5]]: A simplicial homology manifold $K$ of dimension $\geq 5$ 
admits a PL manifold homology resolution iff
\[
\varkappa_4(K) = 0 \;.
\]
The resolutions themselves are classified by maps to $H(\Theta,3)$ [8].
\bigskip

\noindent
{\bf 2.3} Using this technology, the question motivating this note can be formulated as
the problem of understanding the map
\[
B(\tp \La 7 \Ra)  = B(\pl \La 7 \Ra) \to B\comb \;.
\]
Its fiber can be decomposed as
\[   
\pl/\tp \La 7 \Ra \to \comb/\pl \La 7 \Ra \to \comb/\pl = H(\Theta,3) \;;
\]
the fibration
\[
\pl/\tp \La 7 \Ra \to \tp/\tp \La 7 \Ra = H(\Z_2,1) \times H(\Z,3) \to \tp/\pl = H(\Z_2,3)
\]
shows that $\pl/\tp\La 7 \Ra$ is a three-stage Postnikov system, with homotopy
group $\Z$ in degree three and $\Z_2$ in degrees one and two. The group $\pi_*(\comb/\pl \La 7 \Ra)$ is 
therefore $\Z_2$ in degree one and zero in degree two, while in degree three there is an exact sequence
\[
0 \to \Z \to \pi_3(\comb/\tp \La 7 \Ra) \to \Theta \to \Z_2 \to 0 \;.
\]
The map on the right can be identified with Rokhlin's invariant
\[
\Sigma \mapsto \rho(\Sigma) := \frac{1}{8} \; {\rm signature}\; (W) \; {\rm  modulo} \;
(2) : \Theta \to \Z_2
\]
of a homology three-sphere $\Sigma$ (where $W$ is a spin four-manifold with $\partial W =
\Sigma$). Let $\Theta_0$ be the kernel of $\rho$. \medskip

\noindent
{\bf Corollary:} When $K$ is a smoothable ($\varkappa(K) = 0$) simplicial homology string manifold 
of dimension $\geq 5$, PL manifold structures on a homology resolution are classified by elements of
$H^3(|K|,\wTheta)$; where
\[
\wTheta := \Z \oplus \Theta_0 \cong \pi_3(\comb/\tp \La 7 \Ra) \;.
\] 
{\bf Proof:} The exact sequence above splits, because the infinite cyclic group maps isomorphically to
the three-dimensional homotopy group of $\tp$. This suggests, among other things, the existence of 
a combinatorial formula [2] for its Pontrjagin class. When that class vanishes, lifts of the classifying
map from $|K|$ to $B\comb$ to a map from a homology resolution $X$ are classified by maps to $\comb/\tp 
\La 7 \Ra$. $\Box$ \bigskip

\noindent
{\bf 2.4} The inclusion of the fiber in 
\[
H(\Theta,3) \to B\pl \to B\comb
\]
is trivial at odd primes: $B\pl \cong B\otimes$ [30], but $K$-theory is blind to 
Eilenberg - MacLane spaces $H(A,n)$ for $n > 2$. At the prime two, there are still open
questions. In particular, it is not known if $\rho$ splits: this is equivalent to the 
conjecture that all topological manifolds of dim $>$ 4 are simplicial complexes. \medskip

\noindent
Freed [10] identifies the classifying space of the Picard category of $\Z_2$-graded complex lines
as a two-stage Postnikov system. Its associated infinite-loop spectrum
\[
\xymatrix{
\Ell \ar[r] &  H\Z_2 \ar[r]^{\beta Sq^2} & \Sigma^3 H\Z }
\]
is the truncation to positive dimensions
\[
(\Sigma^2 \tI)_{<0} \to \Sigma^2 \tI \to \Ell
\]
of the double suspension of the Anderson dual [17 appendix B, 18] of the sphere spectrum, which is characterized
by a short exact sequence
\[
0 \to {\rm Ext}(E_{*-1},\Z) \to \tI^*(E) \to {\rm Hom}(E_*,\Z) \to 0
\]
associated to a spectrum $E$. This same small Postnikov system appears as the base
\[
F/\pl \cong \Sigma^{4*}H\Z \times \Sigma^{4*+2}H\Z_2 \times \Sigma^2 \Ell \; (* > 1)
\]
of the (two-localization) of the infinite-loop space classifying piecewise-linear structures on a
Poincar\'e-duality space [19]. These observations can then be assembled into the diagram

\[
\xymatrix{
{} & \Sigma^2H\Z_2 = \tp/\pl \ar[dr]^\rho & {} & \Sigma^2\tI \ar[d] \\
\Sigma^3H\Theta \ar[ur] \ar[r]^= & \comb/\pl \ar[r]^\kappa & F/\pl \ar@{.>}[ur] \ar[r] & \Sigma^2\Ell
\;.}
\]

\noindent
On the other hand, the Rokhlin homomorphism defines an exact sequence
\[
\dots \to \tI^2(H\Theta_0) \to \tI^1(H\Z_2) \to \tI^1(H\Theta) \to \dots \;;
\]
but $\tI^2(H\Theta_0) = 0$, so (the nonzero invariant defined by) $\rho \in \tI^1(H\Z_2)$ maps
to (the nonzero invariant defined by) $\kappa \in \tI^1(H\Theta)$: which defines a homomorphism 
from $B\Theta$ to $\Ell$, and can thus be interpreted as a topological field theory mapping 
the category with one object, and 3D homology cylinders as morphisms, to the category of $\Z_2$-graded 
complex lines. \medskip

\noindent
This can be regarded as a lift of Rokhlin's invariant, regarded as a topological field theory 
taking values in the Picard category of real lines. It suggests the interest of super-Chern-Simons
theories [11 \S 9] defined on simplicial homology spin manifolds. \bigskip

\section{\"Uber die Hypothesen, welche zu Grunde der Physik liegen}

\noindent
Following [6 \S VII], it is tempting to interpret the elements of $\Z$ in $\wTheta$ as topological twists 
`in the large' (or at infinity), and elements of $\Theta_0$ as twists `in the small'. Physics has
a tradition of concern (cf eg Wheeler) with the possibility that the microstructure of the Universe
might in some way be nonEuclidean. This seems legitimate: experiments in {\bf very} high-energy physics probe 
the topology of space-time at very short distance, and it is conceivable that at very fine scales physical 
space might be described by some kind of quantum foam model [23], possibly involving ensembles with varying
topology. \medskip
                                  
\noindent 
These ideas have a big literature, but interested researchers seem not to be very aware 
of the long history of interest in analogous questions among topologists. In particular, the
extended homology cobordism group $\wTheta$ seems to capture rather precisely the
idea that the space-time `bubbles' in which very-short-distance interactions occur - I'm thinking
of the way Feynman diagrams are often drawn - might be bounded by non-standard
spheres.\medskip

\noindent
If physics starts by hypothesizing the existence of a simplicial homology manifold 
structure\begin{footnote}{conceptually similar to a causal structure}\end{footnote} on some (say, 
ten or eleven-dimensional) space-time $K$, then the vanishing of $w_2$ decides the existence of a 
spin structure and the vanishing of $\varkappa(K)$ decides the existence of a PL resolution. 
When $p_1/2 = 0$, resolutions of $|K|$ admit string structures; this is quite 
like the standard situation. However, the possible {\bf twists} of that string resolution
(related to $B$-fields [26], Vafa's discrete torsion [2 \S 1.6] and perhaps more generally to $D$-brane charges 
[5 \S 4.4]) lie in $H^3(|K|,\wTheta)$, which is much bigger than the usual group of gerbes over 
$|K|$. There may even be `experimental' evidence for the physical relevance of twisting by homology 
three-spheres, in that deep results about the structure of such manifolds are derived by scattering 
Yang-Mills bosons off them: ie, from Donaldson theory [9]. 
\medskip

\noindent
As for exotic homology manifolds [4], hypotheses non fingo: in
part because they, unlike simplicial homology manifolds, seem to lack the clocks and 
measuring rods that play the role of rulers and compasses in classical geometry. This is probably
a lack of imagination on my part; a deeper concern is that major questions in 4D geometric topology 
remain open. Here I want only to make the point that simplicial homology manifolds are understood 
well enough to test against the models of contemporary physics.

\bibliographystyle{amsplain}

\end{document}